\documentclass{article}
\usepackage{amssymb}
\usepackage{amsmath}
\usepackage{amsfonts}

\newtheorem{theorem}{Theorem}

\newtheorem{corollary}[theorem]{Corollary}

\newtheorem{definition}[theorem]{Definition}
\newtheorem{example}[theorem]{Example}

\newtheorem{lemma}[theorem]{Lemma}

\providecommand{\al}{\alpha}
\providecommand{\be}{\beta}

\begin{document}


\title{Counting The Generator Matrices of $\mathbb{Z}_{2}\mathbb{Z}_{8}$-Codes}
\maketitle
\begin{center}
 Irfan Siap $^{a}$ isiap@yildiz.edu.tr \\
 Ismail Aydogdu $^{a}$ iaydogdu@yildiz.edu.tr\\
$^{a}$ Yildiz Technical University, Faculty of Arts and
Science, Department of Mathematics, Istanbul-Turkey
\end{center}

\begin{abstract}
\noindent In this paper, we count the number of matrices whose rows generate different $\mathbb{Z}_2\mathbb{Z}_8$ additive codes. This is a natural generalization of the well known Gaussian numbers that count the number of matrices whose rows generate vector spaces with particular dimension over finite fields. Due to this similarity we name this numbers as Mixed Generalized Gaussian Numbers (MGN). The MGN formula by specialization leads to the well known formula for the number of binary codes and the number of codes over $\mathbb{Z}_8,$ and for additive $\mathbb{Z}_2\mathbb{Z}_4$ codes. Also, we conclude by some properties and examples of the MGN numbers that provide a good source for new number sequences that are not listed in The On-Line Encyclopedia of Integer Sequences.
\end{abstract}

\begin{quotation}
\end{quotation}

\section{Introduction}

Let $\mathbb{Z}_{m}$   be  the ring of integers modulo $m.$  $\mathbb{Z}_{m}^{n}$ will denote the the set of cartesian product of $n$ copies of $\mathbb{Z}_{m}.$ Any nonempty subset $C$ of $\mathbb{Z}_{m}^{n}$ is called a  code and a subgroup of a $\mathbb{Z}_{m}^{n}$ is called a linear code of length $n.$  For the special cases  $n=2$ and $n=4,$ the codes are called binary and quaternary codes respectively. Most of the work and applications in digital communications is done on binary linear codes. However, due to the relations between the algebraic structures via some special maps which are referred to as Gray maps,  the images of codes over non binary rings  provide structural binary codes. An important such work  is introduced by Hammons et al. \cite{Hammons}  and since then, the study on codes over various rings has been of quite interest on Algebraic Coding Theory. One of such a successful attempt  is the study of codes which are group isomorphic to additive subgroups of the group  $ \mathbb{Z}_{2}^\alpha \times \mathbb{Z}_{4}^\beta $ where $\al $ and $\be $ positive integers. It is clear that in $ \mathbb{Z}_{2}^\alpha \times \mathbb{Z}_{4}^\beta $ if $\be$ does not exist then the subgroups give linear binary codes, or if $\al$ does not exist then the subgroups give linear quaternary codes. So this is a generalization of the well known families of (binary and quaternary) codes  and are known as additive codes. Additive codes were originally defined by Delsarte in 1973 in the context of association schemes \cite{Delsarte 1,Delsarte 2}. Puyol at el. in \cite{Pujol},  translation invariant
propelinear codes are introduced and  these
 codes shown to be  isomorphic to the subgroups of $%
\mathbb{Z}
_{2}^{\alpha }\times
\mathbb{Z}
_{4}^{\beta }\times
\mathbb{Q}
_{8}^{\sigma }$ where $%
\mathbb{Q}
_{8}$ is the nonabelian quaternion group with eight elements. In the binary Hamming scheme case,  the additive codes are exactly the abelian translation invariant propelinear codes. Therefore, the
only structures for the abelian groups are  $%
\mathbb{Z}
_{2}^{\alpha }\times
\mathbb{Z}
_{4}^{\beta },$ where $\alpha +2\beta =n~\cite{Pujol}.$  Also, codes defined over two different alphabets which are binary and ternary fields is studied by Brouwer at el. in \cite{mixed}.

The basic and introductory concepts on  $%
\mathbb{Z}
_{2}%
\mathbb{Z}
_{4}-$codes  are presented in \cite{eskiZ2Z4, deskiZ2Z4,Borges} where the readers can further refer to.  An additive code $C$ over $\mathbb{Z}_{2}\times \mathbb{Z}_{4}$ which is a subgroup of $\mathbb{Z}_{2}^\al\times \mathbb{Z}_{4}^\be$ is group isomorphic to $\mathbb{Z}_{2}^{k_0}\mathbb{Z}_{2}^{k_1}\times \mathbb{Z}_{4}^{k_2}.$ Here, $k_0$ represents the number of generators of the subgroup $C$ of order 2 that are contributed through the binary ($\mathbb{Z}_2$) part, $k_1$ represents the number of generators of the subgroup $C$ of order 4 that are contributed through the quaternary ($\mathbb{Z}_4$) part and $k_2$ represents the number of generators of the subgroup $C$ of order 2 that are contributed through the quaternary ($\mathbb{Z}_4$) part. Thus, this leads to the following fact that is proved in \cite{Borges}: an additive $\mathbb{Z}_{2} \mathbb{Z}_{4}$ code of type $(\alpha ,\beta ;k_{0},k_{1},k_{2})$ is equivalent to an additive code generated by the following matrix (\cite{Borges})

\begin{equation} \label{MZ2Z4}
G=%
\left[\begin{array}{cc|ccc}
I_{k_{0}} & \bar{A}_{01} & 0 & 0 & 2T_{02} \\
0 & S_{1} & I_{k_{1}} & A_{01} & A_{02} \\
0 & 0 & 0 & 2I_{k_{2}} & 2A_{12}%
\end{array}\right]%
,
\end{equation}%
(P.S. The artificial vertical  line only helps to distinguish between the binary and the quaternary parts).

Similar to the discussions above, if $C$ is a $%
\mathbb{Z}
_{2}%
\mathbb{Z}
_{8}-$ additive code of type $(\alpha ,\beta ;k_{0},k_{1},k_{2},k_{3}), $  then in \cite{amis} it is proven that   $C$ is equivalent to a code generated by the following matrix (\cite{amis})

\begin{equation}\label{MZ2Z8}
G=%
\left[\begin{array}{cc|cccc}
I_{k_{0}} & \bar{A}_{01} & 0 & 0 & 0 & 4T_{03} \\
0 & S_{1} & I_{k_{1}} & A_{01} & A_{02} & A_{03} \\
0 & S_{2} & 0 & 2I_{k_{2}} & 2A_{12} & 2A_{13} \\
0 & 0 & 0 & 0 & 4I_{k_{3}} & 4A_{23}%
\end{array}\right].
\end{equation}

Here $k_0$ represents the number of order 2 generators that are contributed through the binary part, and respectively, $k_1,$ $k_2$ and $k_3$ represent the number of order 8, 4 and 2 generators  2 that are contributed through the $\mathbb{Z}_8$ part. Note that the order 2 elements from the $\mathbb{Z}_8$ part have the first $\alpha$ components all zero. This remark will play a crucial role in the main counting theorem in the next section.

We present some facts regarding the duality of this codes which is introduced in \cite{amis}.The inner product of two vectors $u,v\in
\mathbb{Z}
_{2 }^\al\times
\mathbb{Z}
_{8}^{\beta }$ \bigskip is defined as;

\begin{equation*}
\langle u,v \rangle =4\left( \sum_{i=1}^{\alpha }u_{i}v_{i}\right) +\sum_{j=\alpha
+1}^{\alpha +\beta }u_{j}v_{j}\in
\mathbb{Z}
_{8}.
\end{equation*}
The additive dual code of $C,$ denoted by $C^{\perp
},~$ is then defined as

\begin{equation*}
C^{\perp }=\left\{ v\in
\mathbb{Z}
_{2}^{\alpha }\times
\mathbb{Z}
_{8}^{\beta }|\langle u,v\rangle =0~\text{for all }u\in C\right\} .
\end{equation*}

It is easy to check that $C^{\perp }$ is a subgroup of $%
\mathbb{Z}
_{2}^{\alpha }\times
\mathbb{Z}
_{8}^{\beta },$ so $C^{\perp }$ is a $%
\mathbb{Z}
_{2}%
\mathbb{Z}
_{8}-$additive code too.

Let $C$ be a $%
\mathbb{Z}
_{2}%
\mathbb{Z}
_{8}-$additive code of type $(\alpha ,\beta ;k_{0},k_{1},k_{2},k_{3}) $ with canonical generator matrix \ref{MZ2Z8}. Then, the parity-check
matrix of $C$  which is the generator matrix of its dual is
\begin{equation*}
\begin{bmatrix}
-\bar{A}_{01}^{t} & I_{\alpha -k_{0}} &  & -2S_{2}^{t} & 0 & 0 \\
-T_{03}^{t} & 0 & P & -A_{13}^{t}+A_{23}^{t}A_{12}^{t} & -A_{23}^{t} &
I_{\beta -k_{1}-k_{2}-k_{3}} \\
0 & 0 &  & -2A_{12}^{t} & 2I_{k_{3}} & 0 \\
0 & 0 &  & 4I_{k_{2}} & 0 & 0%
\end{bmatrix}%
,
\end{equation*}

where $P=%
\begin{bmatrix}
-4S_{1}^{t}+2S_{2}^{t}A_{01}^{t} \\
-A_{03}^{t}+A_{13}^{t}A_{01}^{t}+A_{23}^{t}A_{02}^{t}-A_{23}^{t}A_{12}^{t}A_{01}^{t}
\\
-2A_{02}^{t}+2A_{12}^{t}A_{01}^{t} \\
-4A_{01}^{t}%
\end{bmatrix}%
.$

If $C$ is an $%
\mathbb{Z}
_{2}%
\mathbb{Z}
_{8}-$ additive code of type $(\alpha ,\beta ;k_{0},k_{1},k_{2},k_{3}), $  then $C^\perp$ is an $%
\mathbb{Z}
_{2}%
\mathbb{Z}
_{8}-$ additive code of type $(\alpha ,\beta ; \alpha -k_{0},\beta -k_{1}-k_2-k_3,k_{3},k_{2}).$

\section{ Mixed Generalized Gaussian Numbers }

In this section, we present the main theorem of this paper that gives a direct computation of the number of matrices that generate different (not necessarily equivalent) additive codes. First, we present a very moderate example in order to illustrate the problem. Even in this example, getting the exact number may not be easy. As the size of the matrix gets larger the difficulty of counting these matrices becomes a very difficult problem. After stating and proving the main theorem we revisit this example and solve the problem directly.

 As mentioned in introduction the counting problem is originated from the study of the number of the subspaces generated by the rows of matrices over finite fields. Recently, there has been some generalizations of these concept on the number of generating matrices of particular types over the ring $\mathbb{Z}_m$ \cite{S}, over Galois rings \cite{S} and over rings $F_p+uF_p.$ The main Theorem \ref{main} presents a further generalization to the work done in \cite{S} in a different direction.

\begin{example}\label{ex22110}
Let $C$ be a $%
\mathbb{Z}
_{2}%
\mathbb{Z}
_{8}-$ additive code of type $(2,2;1,1,1,0)$ then all possible matrices are 36 matrices that generate different codes.
Here $\alpha=2, \beta =2, k_0=k_1=k_2=1,$ and $k_3=0.$

\begin{enumerate}
\item
$\begin{bmatrix}
1 & g_{21} & 0 & 0  \\
0 & g_{22} &1 & g_{24} \\
0 & g_{32} & 0 & 2%
\end{bmatrix}
$  here we have 16 possible matrices.

\item
$\begin{bmatrix}
1   & g_{21} & 0   & 0  \\
0   & g_{22} &0    & 1 \\
0   & g_{32}      & 2     & 0%
\end{bmatrix}
$  here we have 8 possible matrices.

\item
$\begin{bmatrix}
0         & 1 & 0        & 0  \\
g_{21}   & 0  &1    & g_{24}\\
g_{31}         & 0      & 0     & 2%
\end{bmatrix}
$  here we have 8 possible matrices.

\item
$\begin{bmatrix}
0         & 1 & 0        & 0  \\
g_{21}   & 0  &0    & 1 \\
g_{31}       & 0      & 2     & 0%
\end{bmatrix}
$  here we have 4 possible matrices.

\end{enumerate}

where all unknown above are either $0$ or $1.$ So altogether we have $36$ generating matrices.

\end{example}

\begin{theorem}\label{main}

Let

$$N_1(\alpha,\beta;k_0,k_1,k_2,k_3)=\prod_{i=0}^{k_0-1}(2^\alpha-2^i)2^\beta,$$ (number of choices to form an element of order 2 that is contributed through the binary part from all space);

$$N_2(\alpha,\beta;k_0,k_1,k_2,k_3)=\prod_{i=0}^{k_1-1}(8^{\beta}-4^{\beta }\cdot 2^i) 2^\alpha,$$ (number of choices to form an element of order 8  that is contributed through the $\mathbb{Z}_8$ part from all space);

$$N_3(\alpha,\beta;k_0,k_1,k_2,k_3)=\prod_{i=0}^{k_2-1}(4^\beta-2^{\beta +k_{1}+i}) 2^\alpha ,$$ (number of choices to form an element of order 4  that is contributed through the $\mathbb{Z}_8$ part from all space); and

$$N_4(\alpha,\beta;k_0,k_1,k_2,k_3)=\prod_{i=0}^{k_3-1}(2^{\beta}-2^{{k_2} +{k_1}+i}) ,$$ (number of choices to form an element of order 2  that is contributed through the $\mathbb{Z}_8$ part from all space).

Let

$$D_1(\alpha,\beta;k_0,k_1,k_2,k_3)=\prod_{i=0}^{k_0-1}(2^{k_0+k_1+k_2+k_3}-2^{k_1+k_2+k_3+i}),$$ (number of choices to form an element of order 2 that is contributed through the binary part within the code of the given type);

$$D_2(\alpha,\beta;k_0,k_1,k_2,k_3)=\prod_{i=0}^{k_1-1}(8^{k_1}-4^{k_1}2^i)2^{k_0+2k_2+k_3},$$ (number of choices to form an element of order 8  that is contributed through the $\mathbb{Z}_8$ part  within the code of the given type);

$$D_3(\alpha,\beta;k_0,k_1,k_2,k_3)=\prod_{i=0}^{k_2-1}(4^{k_2}-2^{k_2+i})2^{k_0+2k_1+k_3},$$ (number of choices to form an element of order 4  that is contributed through the $\mathbb{Z}_8$ part  within the code of the given type); and

$$D_4(\alpha,\beta;k_0,k_1,k_2,k_3)=\prod_{i=0}^{k_3-1}(2^{k_1+k_2+k_3}-2^{k_1+k_2+i}),$$ (number of choices to form an element of order 2  that is contributed through the $\mathbb{Z}_8$ part  within the code of the given type).

The number of $\mathbb{Z}_2 \mathbb{Z}_8$ additive codes of the type $(k_0,k_1,k_2,k_3)$ is equal to

\begin{equation} \label{M}  N_{2\times 8}(\alpha,\beta;k_0,k_1,k_2,k_3)=\frac{N_1N_2N_3N_4}{D_1D_2D_3D_4}.\end{equation}

(P.S. When it is clear from the context we will drop the notation $(\alpha,\beta;k_0,k_1,k_2,k_3)$ that explains each computation. Further if any of the parameters $k_i; (0\leq i\leq 3)$ is zero then the $N_i=1$ and $D_i=1.$)

\end{theorem}

{\bf Proof:} In order to prove this theorem we count ordered generators for the group (code) of the type $(k_0,k_1,k_2,k_3).$ First we count them by choosing them from the all space $\mathbb{Z}_2^\al\times \mathbb{Z}_8^\be$ which gives say $A$ and given a group (code) of type   $(k_0,k_1,k_2,k_3)$ then we choose them in within this group. So, if the number of groups of  the type $(k_0,k_1,k_2,k_3)$ is $N_{2\times 8},$ then $N_{2\times 8}=A/B.$ First, we compute $A:$ In $\mathbb{Z}_2^\al\mathbb{Z}_8^\be$ we can choose an element of order $2$  that is contributed through the binary part in $(2^\al-1)\cdot 2^\be.$ Next, the second element with the same property can be choose in   $(2^\al-2)\cdot 2^\be$ ways,  inductively the last element can be chosen in $(2^\al-2^{k_0-1})\cdot 2^\be.$ So, in total $k_0$ elements of order $2$  that  contribute through the binary part in all space can be chose in $N_{1}$ ways. Next, there are $(8^n-4^n)\cdot 2^\al$ ways to pick an element of order $8$ contributed through the $\mathbb{Z}_8$ part. The second such element can be chosen in $(8^n-4^n\cdot 4)\cdot 2^\al$ ways excluding the linear combinations of the first chosen element of order $8.$ Inductively, we have $N_2$ choices for such elements. Next, to choose elements of order $4$ in the all space, first there are $(4^n-2^n)\cdot 2^{k_1}2^\al$ to pick elements of order $4$ that are contributed through the $\mathbb{Z}_8$ part. Here, the elements of order $4$ that are formed from the $k_1$ elements of order $8$ by taking their $2$ multiples need to be considered. Then, similarly as discussed above we have $N_3$ elements of order $4$ in the all space to be chosen. Next, to choose elements that are of order $2$ and solely contributed through $\mathbb{Z}_8$ part. This imposes that the first $\al$ entries of such elements to be all zero. Thus, in order to pick such an element first we subtract order $2$ elements that are obtained through already chosen $k_1$ and $k_2$ elements for the $\mathbb{Z}_8$ part. So, we have $(2^{k_0+k_1+k_2+k_3}-2^{k_0+k_1+k_2})$ choices. The next choice comes from  $(2^{k_0+k_1+k_2+k_3}-2\cdot 2^{k_0+k_1+k_2})$ and inductively we reach $N_4.$ hence, $A=N_1N_2N_3N_4.$ Now, we repeat the same process within the group of the type $(k_0,k_1,k_2,k_3).$ In order to pick an element of order $2$ that   contributes through the binary part, we subtract all order $2$ elements within the group from the ones that are not coming through the binary part, i.e $(2^{k_0+2k_1+k_2+k_3}-2\cdot 2^{k_0+2k_1+k_2})$  and inductively we reach at $D_1.$ Next, to choose an element within the group of order $8,$ we have  $(8^{k_1}\cdot 2^{k_0+2k_2+k_3}-4^{k_1}\cdot 2^{k_0+2k_2++k_3})$ choices. Inductively, we have $D_2$ choices altogether. Next, to choose an element within the group of order $4,$ we have  $(4^{k_2}-2^{k_2})\cdot 2^{k_0+2k_1+k_3}$ choices. Inductively, we have $D_3$ choices altogether. Finally, to choose an element within the group of order $2,$ we have  $(2^{k_1+k_2+k_3}-2^{k_1+k_2})$ choices.  Inductively, we have $D_4$ choices altogether. Therefore, we have the result. $\qquad \square$

\begin{example}(An application of the main theorem) In Example \ref{ex22110}, by explicitly working out the cases,  we computed the $N_{2\times 8}(2,2;1,1,1,0)$
 which is alternatively given in a direct way by Theorem \ref{M}:

 $N_1(2,2;1,1,1,0)=12,$  $N_2=192,$ $N_3=32,$ $N_4$ does not exist so we skip this term. $D_1=4,$  $D_2=32,$  $D_3=16,$ and we disregard $D_4.$ Hence,
 $$N_{2\times 8}(2,2;1,1,1,0)=\frac{N_1N_2N_3}{D_1D_2D_3}=\frac{73728}{2048}=36.$$

\end{example}

The next corollary shows that the number of distinct linear codes over $\mathbb{Z}_8$ (\cite{Sa}) can be obtained via the main Theorem \ref{main}.

\begin{corollary}\label{coroZ8}  Let $N_8(n;k_1,k_2,k_3)$ be the number of distinct linear codes of type $(k_1,k_2,k_3)$ over $\mathbb{Z}_8.$  If, $r=n-k_1+k_2+k_3,$ then,

$$N_{8}(n,k_1,k_2,k_3)=  2^{-r}\cdot N_{2\times 8}(1,n;1,k_1,k_2).$$
\end{corollary}
In the following we present some applications of Corollary \ref{coroZ8}.
\begin{example}
$$ N_{8}(4;2,1,1)= N_{2\times 8}(1,4;1,2,1,1)=420.$$
$$ N_{8}(5;3,0,0)=  2^{-2}\cdot N_{2\times 8}(1,5;1,3,0,0)=2539520/4=634880.$$
$$ N_{8}(6;2,0,1)=  2^{-3}\cdot N_{2\times 8}(1,6;1,2,0,1)=159989760.$$
\end{example}

\section{Some Connections}

In this section, we explore some connections of this formula to other related topics.

\subsection{$q$-binomial and multinomial coefficients}

In this subsection we introduce some basic well known facts regarding the $q$-binomials. The main goal is to explicitly write down an identity that relates the $q$-binomial coefficients with mixed Gaussian numbers.

Let $n$  and $q$ be two positive integers.

$$[n]_q=1+q+q^2+\cdots +q^{n-1}=\frac{q^n-1}{q-1}.$$

The $q$ -factorial  is defined as

$$[n]_q!=[n]_q\cdot [n-1]_q\cdots [2]_q\cdot [1]_q.  $$

\begin{definition}\cite{bona}($(n,k)$-Gaussian coefficient or $q$-binomial)  Let $n,k$  and $q$ be   non negative integers  such that $k\leq n.$ Then,
\begin{equation} \label{q} \left[\begin{array}{c}
n\\
k
\end{array}
\right]_q=\frac{[n]_q!}{[k]_q![n-k]_q!}.
\end{equation}
Here, $[0]_q!=1$  and  if $k=0$ then the $q$ binomial coefficient is equal to $1.$
\end{definition}

The formula (\ref{q}) has many interpretations. An algebraic interpretation is that it  gives the number of vector subspaces of dimension $k$ in $F_q^n$ where $F_q$ is a finite field of order $q.$  This is also equivalent to the number of matrices on the row equivalent form of size $k\times n$ and rank $k.$  Also, combinatorially, (\ref{q}) can also be interpreted as a polynomial in $q$  where the coefficient  $q^k$ counts the number of distinct partitions of  $k$ elements which fit inside a  rectangle of size $k \times (n-k)$. A nice survey and a different interpretation of $q$ binomial coefficients via tiling can be found in \cite{tilings}. Further, some relations between $q$ binomial and classical binomial coefficients, and some special polynomials with their applications to distributions are exposed in \cite{alger}.

The very first two well known properties of $q$ binomial numbers are listed below:

\begin{equation}\label{PropQ} \left[\begin{array}{c}
n \\
k
\end{array}
\right]_q = \left[\begin{array}{c}
n \\
n-k
\end{array}
\right]_q  \text{ and } \left[\begin{array}{c}
n \\
k
\end{array}
\right]_q=q^k \left[\begin{array}{c}
n - 1\\
k
\end{array}
\right]_q+ \left[\begin{array}{c}
n - 1\\
k-1
\end{array}
\right]_q.
\end{equation}

A further and natural generalization of $q$ binomial coefficients is $q$ multinomial coefficients. Similar to binomial coefficients $q$ multinomial coefficients have found applications on distribution and hence statistics and physics \cite{ph}. An algebraic expression of these numbers is given in \cite{Vmulti} and \cite{flag} as follows: Let $V=F_q^n$ be a vector space of dimension $n$ over a finite field $F_q$ with $q$ elements. Let $V_i$ be subspaces of $V$ each of dimension $n-\sum_{j=1}^ik_i$ where $1\leq i \leq m$ and $V_1\subset V_2 \subset \cdots \subset V_m\subset V.$  Such a chain of subvector spaces is referred to as a flag of subspaces of length $m$  of $V.$ Then, the number of a flag of subspaces of length $m$ specified above is given by the $q$ multinomial coefficient:

\begin{equation}\label{Vmulti} \left[\begin{array}{c}
n \\
k_1 ,k_2, \ldots  ,k_m
\end{array}
\right]_2 = \left[\begin{array}{c}
n \\
k_1
\end{array}
\right]_2 \cdot \left[\begin{array}{c}
n -k_1 \\
k_2
\end{array}
\right]_2  \cdots \cdot \left[\begin{array}{c}
n - \sum_{i=1}^{m-1} k_i\\
k_{m}
\end{array}
\right]_2
\end{equation}
where $n=\sum_{i=1}^m k_i.$

In the sequel we relate the $q$ binomial and multinomial coefficients with MGN. First, the following corollary shows that the number of distinct binary linear codes over $\mathbb{Z}_2$ (\cite{sloane}) of length $n$ can be obtained via the main Theorem \ref{main}, by simple observation we get the formula for the $2$-binomial coefficients (Gaussian Numbers over $\mathbb{Z}_2$).

\begin{corollary}\label{coroZ2}
$\left[\begin{array} {cc} n\\ k \end{array}\right]_2=N_{2\times 8}(n,1;k,0,0,1).$
\end{corollary}

Now, we relate the main formula (\ref{main}) with $2$-binomial coefficients. To  accomplish this task we express  the formula for Mixed Gaussian Numbers by standard $2$ binomial coefficients. First, we observe the following:

\begin{align*}\prod_{i=0}^{t-1}(2^r-2^i)&=(2^r-1)(2^r-2)(2^r-2^2)\cdots (2^r-2^{t-2})(2^r-2^{t-1})\\
                                        &=2^{1+2+\cdots +t-1}(2^r-1)(2^{r-1}-1)(2^{r-2}-1)\cdots (2^{r-(t-1)}-1)\\
                                        &=2^{\binom{t}{2}}\frac{[r]_2!}{[r-t]_2!}.
\end{align*}

Next, we express all elements in the main formulae by $2$ binomial coefficients:

$$N_1=\prod_{i=0}^{k_0-1}(2^\alpha-2^i)2^\beta=2^{k_0\beta}2^{\binom{k_0}{2}}\frac{[\alpha]_2!}{[\alpha -k_0]_2!},  $$

$$N_2=\prod_{i=0}^{k_1-1}(8^{\beta}-4^{\beta }\cdot 2^i) 2^\alpha=2^{2\beta k_1+\alpha k_1 }  \prod_{i=0}^{k_1-1}(2^{\beta}- 2^i)=2^{2\beta k_1+\alpha k_1+\binom{k_1}{2}}\frac{[\beta]_2!}{[\beta -k_1]_2!},  $$

$$N_3=\prod_{i=0}^{k_2-1}(4^\beta-2^{\beta +k_{1}+i}) 2^\alpha =2^{(\beta+ k_1)k_2+\alpha k_2 }  \prod_{i=0}^{k_2-1}(2^{\beta-k_1}- 2^i)=2^{(\beta + k_1)k_2+\alpha k_2+\binom{k_2}{2}}\frac{[\beta-k_1]_2!}{[\beta - k_1-k_2]_2!},  $$

$$N_4=\prod_{i=0}^{k_3-1}(2^{\beta}-2^{{k_2} +{k_1}+i})=2^{( k_1+k_2)k_3 }  \prod_{i=0}^{k_3-1}(2^{\beta-k_1-k_2}- 2^i)=2^{( k_1+k_2)k_3 +\binom{k_3}{2}}\frac{[\beta-k_1-k_2]_2!}{[\beta - k_1-k_2-k_3]_2!},  $$

$$D_1=\prod_{i=0}^{k_0-1}(2^{k_0+k_1+k_2+k_3}-2^{k_1+k_2+k_3+i})=2^{(k_1+k_2+k_3)k_0 }  \prod_{i=0}^{k_0-1}(2^{k_0}- 2^i)=2^{(k_1+k_2+k_3)k_0 +\binom{k_0}{2}}{[k_0]_2!},  $$

$$D_2=\prod_{i=0}^{k_1-1}(8^{k_1}-4^{k_1}2^i)2^{k_0+2k_2+k_3}=2^{(k_0+2k_1+2k_2+k_3)k_1 }  \prod_{i=0}^{k_1-1}(2^{k_1}- 2^i)=2^{(2k_1+k_0+2k_2+k_3)k_1 +\binom{k_1}{2}}{[k_1]_2!},  $$

$$D_3=\prod_{i=0}^{k_2-1}(4^{k_2}-2^{k_2+i})2^{k_0+2k_1+k_3}=2^{(k_0+2k_1+k_2+k_3)k_2 }  \prod_{i=0}^{k_2-1}(2^{k_2}- 2^i)=2^{(2k_1+k_0+k_2+k_3)k_2 +\binom{k_2}{2}}{[k_2]_2!} , $$

and

$$D_4=\prod_{i=0}^{k_3-1}(2^{k_1+k_2+k_3}-2^{k_1+k_2+i})=2^{(k_1+k_2)k_3}  \prod_{i=0}^{k_3-1}(2^{k_3}- 2^i)=2^{(k_1+k_2)k_3 +\binom{k_3}{2}}{[k_3]_2!}.  $$

Finally,

\begin{equation*}
N=2^\delta \frac{[\alpha]_2!}{[k_0]_2! [\alpha -k_0]_2!}\frac{[\beta]_2!}{[k_1]_2![\beta -k_1]_2!} \frac{[\beta -k_1]_2!}{[k_2]_2! [\beta -k_1-k_2]_2!}   \frac{[\beta-k_1-k_2]_2!}{[k_3]_2![\beta - k_1-k_2-k_3]_2!}
\end{equation*}

i.e.

\begin{equation}\label{eqConn} N=2^\delta \left[\begin{array}{c}
\alpha\\
k_0
\end{array}
\right]_2 \cdot \left[\begin{array}{c}
\beta \\
k_1
\end{array}
\right]_2 \cdot \left[\begin{array}{c}
\beta -k_1 \\
k_2
\end{array}
\right]_2 \cdot \left[\begin{array}{c}
\beta - k_1 - k_2\\
k_3
\end{array}
\right]_2
\end{equation}

and further,  we have

\begin{equation}\label{eqConn} N=2^\delta \left[\begin{array}{c}
\alpha\\
k_0
\end{array}
\right]_2 \cdot \left[\begin{array}{c}
\beta \\
k_1 ,k_2, k_3
\end{array}
\right]_2
\end{equation}

where $\delta = k_0(\beta -l)+k_1(\alpha -k_0+2(\beta-l)+k_3) +k_2((\beta -l) +(\alpha -k_0)$ and $l=k_1+k_2+k_3.$ The  right end side of $N$ is shown to be
 a $2^\delta$ multiple of a special $q$ multinomial.

So, we can state the following result.

\begin{lemma}\label{lemMm}
\begin{equation} \label{Mm}  N_{2\times 8}(\alpha,\beta;k_0,k_1,k_2,k_3)=2^\delta \left[\begin{array}{c}
\alpha\\
k_0
\end{array}
\right]_2 \cdot \left[\begin{array}{c}
\beta \\
k_1 ,k_2, k_3
\end{array}
\right]_2
\end{equation}
where $\delta = k_0(\beta -l)+k_1(\alpha -k_0+2(\beta-l)+k_3) +k_2((\beta -l) +(\alpha -k_0))$ and $l=k_1+k_2+k_3.$
\end{lemma}

Hence the number of distinct codes of a  code having dual code parameters of a code of type $(\alpha ,\beta ;k_{0},k_{1},k_{2},k_{3}), $  is then equal to

\begin{equation} \label{PDual}  N_{2\times 8}(\alpha ,\beta ; \alpha -k_{0},\beta -l,k_{3},k_{2})=2^{\overline \delta} \left[\begin{array}{c}
\alpha\\
\alpha -k_0
\end{array}
\right]_2 \cdot \left[\begin{array}{c}
\beta \\
\beta -l , k_3, k_2
\end{array}
\right]_2
\end{equation}
where $l=\sum_{i=0}^2k_i$  and $\overline \delta =k_1(\alpha -k_0) +(\beta -l)(k_0 +2k_1 +k_2) + k_3(k_1+k_0)$.

The following lemma  states a condition for the number of  codes that equal to the number their duals can be shown by applying the definitions carefully:

\begin{lemma}\label{lemeqDaul} If  $\alpha k_2=k_0(k_2+k_3),$ then
\begin{equation} \label{PDual}  
N_{2\times 8}(\alpha ,\beta ; \alpha -k_{0},\beta -l,k_{3},k_{2})=N_{2\times 8}(\alpha ,\beta ; \alpha -k_{0},\beta -l,k_{3},k_{2}).
\end{equation} 
\end{lemma}

\vspace{2mm}
As special cases to the previous Lemma \label{leneqDual}, we have the following two corollaries:

\begin{corollary} If $C$ is an additive $\mathbb{Z}_{2} \mathbb{Z}_{8}$-code of type $(r,s;k_0,k_1,0,0),$ then the number of such codes is equal to the number of codes with parameters of its dual.
\end{corollary}
\begin{corollary} The number of additive $\mathbb{Z}_{2} \mathbb{Z}_{8}$-codes that are of type $(r,s;k_0,0,k_2,s-k_2)$ and their duals are equal.
\end{corollary}

\subsection{A Connection to $\mathbb{Z}_2 \mathbb{Z}_4$ Additive Codes}

Here, we relate the additive codes over $\mathbb{Z}_{2}\times \mathbb{Z}_{8}$ with additive codes over  $\mathbb{Z}_{2}\times \mathbb{Z}_{4}$ and hence we also obtain a formula for counting the number of the matrices that generate additive codes over $\mathbb{Z}_{2}\times \mathbb{Z}_{4}.$

First we present an explicit example that lists and counts all matrices that generate $\mathbb{Z}_{2} \mathbb{Z}_{4}$ additive codes of type  $(\alpha=2, \beta=2; k_0=k_1=k_2=1).$

\begin{example}\label{ex2211}
 The generators of $\mathbb{Z}_{2} \mathbb{Z}_{4}$ additive codes of type  $(2,2;1,1,1)
$
are

\begin{enumerate}
\item
$\begin{bmatrix}
1 & g_{21} & 0 & 0  \\
0 & g_{22} &1 & g_{24} \\
0 & 0 & 0 & 2%
\end{bmatrix}
$  here we have 8 possible matrices.

\item
$\begin{bmatrix}
1   & g_{21} & 0   & 0  \\
0   & g_{22} &0    & 1 \\
0   & 0      & 2     & 0%
\end{bmatrix}
$  here we have 4 possible matrices.

\item
$\begin{bmatrix}
0         & 1 & 0        & 0  \\
g_{21}   & 0  &1    & g_{24}\\
0        & 0      & 0     & 2%
\end{bmatrix}
$  here we have 4 possible matrices.

\item
$\begin{bmatrix}
0         & 1 & 0        & 0  \\
g_{21}   & 0  &0   & 1 \\
0        & 0      & 2     & 0%
\end{bmatrix}
$  here we have 2 possible matrices.

\end{enumerate}

where all unknown above are either $0$ or $1.$ So, altogether there are 18 such matrices.

\end{example}

 As we see from the above example listing and then counting matrices of special type is a very challenging problem. To avoid such an approach and to achieve this goal first we define an auxiliary map $\phi $ from $\mathbb{Z}_{2}^s\times \mathbb{Z}_{8}^s$ to itself such that $(a_1,\ldots a_r,b_1,\ldots b_2) \in \mathbb{Z}_{2}^s\times \mathbb{Z}_{8}^s$ and
$\phi(a_1,\ldots a_r,b_1,\ldots b_2)=(a_1,\ldots a_r,b_1 \mod 4,\ldots  b_s \mod 4).$ If $C$ is a $\mathbb{Z}_{2}\mathbb{Z}_{8}$ code and $\phi(C) $ is also a  $\mathbb{Z}_{2} \mathbb{Z}_{8}$ additive code. Suppose that $C$ is of type $(\alpha, \beta; k_0,k_1,k_2,k_3)$ which is module isomorphic to a  $\mathbb{Z}_{2}\mathbb{Z}_{4}$ code of type $(\alpha, \beta; k_0,k_1,k_2).$

Let $ N_{2\times 4}(\alpha,\beta;k_0,k_1,k_2)$ denote the number of $\mathbb{Z}_{2}\mathbb{Z}_{4}$ distinct additive codes of type $(\alpha,\beta;k_0,k_1,k_2).$

By making use of the observation and facts mentioned above we easily obtain the following corollary that gives a formula for the number of $\mathbb{Z}_{2}\mathbb{Z}_{4}$ additive codes of the type $(k_0;k_1,k_2)$ by making use of Theorem \ref{main}.

\begin{corollary}\label{coroZ2Z4}  $$ N_{2\times 4}(\alpha,\beta;k_0,k_1,k_2)= N_{2\times 8}(\alpha,\beta;k_0,0,k_2,k_3).$$
\end{corollary}

\begin{example}$$ N_{2\times 4}(3,4;2,1,2)= N_{2\times 8}(3,4;2,0,1,2)=11760.$$
\end{example}

\section{Some Properties And  New Number Sequences}

  The main theorem produced a formula that also enjoys some  properties by its own  such as classical and Gaussian binomials do. Here, we present some properties and also some new number sequences that are not recorded yet in the literature \cite{numbers}.

\begin{lemma}\label{lemD}

Let $r,s,k,l,m,t$ be non negative integers. Also let $m\leq r$ and $s=k+l.$  Then,
    \begin{equation*} N_{2\times 8}(r,s; m,k,l,0)=N_{2\times 8}(r,s; m,l,k,0).\end{equation*}
\end{lemma}

{\bf Proof:}  \begin{align*}N_{2\times 8}(r,s; m,k,l,0)=&\left[\begin{array}{c} r \\  m \end{array}\right]_2\left[\begin{array}{c} s \\  s-k \end{array}\right]_2\left[\begin{array}{c} s-k \\  l \end{array}\right]_2\\
=&\left[\begin{array}{c} r \\  m \end{array}\right]_2\left[\begin{array}{c} s \\ l \end{array}\right]_2\left[\begin{array}{c} s-l \\  k \end{array}\right]_2=N_{2\times 8}(r,s; m,l,k,0)
\end{align*}
since (\ref{PropQ}) and the last binomials at the very right of both lines of equations are equal to one due to $s=k+l.  \square$

\begin{lemma}

Let $r,s\in \mathbb{Z}^+.$
\begin{enumerate}
\item
    $\frac{ N_{2\times 8}(r+1,s; 1,1,1,0)}{N_{2\times 8}(r,s; 1,1,1,0)}=4\frac{(2^{r+1}-1)}{(2^{r}-1)}.$

\item $ N_{2\times 8}(r,s; r,s,0,0)=1,  N_{2\times 8}(r,s; r,0,s,0)=1, \text{ and } N_{2\times 8}(r,s; r,0,0,s)=1.$

\item $N_{2\times 8}(1,r;1,1,1,0)=2^{4r-9}(2^r-2)(2^r-1),$
$$N_{2\times 8}(\alpha+1,r; 1,1,1,0)=4N_{2\times 8}(\alpha,r; 1,1,1,0)+(2^r-1)\cdot (2^{r-1}-1)\cdot 2^{3\alpha+4(r-2)}$$ where $\alpha \geq 1, r \geq 2.$
\item $N_{2\times 8}(\alpha,\beta; \alpha,k_1,k_2,k_3)=N_{2\times 8}(1,\beta; 1,k_1,k_2,k_3)$ for all $\al \geq 1.$

\end{enumerate}

\end{lemma}

{\bf Proof: }

\begin{enumerate}
\item  \begin{align*} N_{2\times 8}(r+1,s; 1,1,1,0)=&\frac{(2^{r+1}-1)2^s\cdot (8^s-4^s)2^{r+1} \cdot (4^s-2^s\cdot 2)2^{r+1}}{(2^{4}-2^3)\cdot (8^1-4^1)2^{3} \cdot (4^1-2^1)2^{3}}\\
&=\frac{(2^{r+1}-1)2^s\cdot (8^s-4^s)2^{r+1} \cdot (4^s-2^s\cdot 2)2^{r+1}}{2^{12}}.\\
\end{align*}

On the other hand,

\begin{align*} N_{2\times 8}(r,s; 1,1,1,0)=&\frac{(2^{r}-1)2^s\cdot (8^s-4^s)2^{r} \cdot (4^s-2^s\cdot 2)2^{r}}{(2^{4}-2^3)\cdot (8^1-4^1)2^{3} \cdot (4^1-2^1)2^{3}}\\
&=\frac{(2^{r}-1)2^s\cdot (8^s-4^s)2^{r} \cdot (4^s-2^s\cdot 2)2^{r}}{2^{12}}.\\
\end{align*}

Hence, we have the result.

\item  By definitions, it is straightforward to see that $N_1=D_1$ and $N_2=D_2,$ hence $N_{2\times 8}(r,s; r,s,0,0)=\frac{N_1N_2}{D_1D_2}=1.$ The other two identities follow in a similar way.
    \item By definitions it follows.
\item By definitions it follows. $\qquad \square$
 \end{enumerate}

Besides many sequences that we run into in  this research we would like to mention a few of them.  $N_{2\times 8}(1,k;1,1,1,1)$ where $k\geq 3$  and the sequence with its first three entries is $\{42,10080,1666560,239984640, \ldots\}.$ This also a new sequence which is not listed in \cite{numbers}. We present some new sequences that are not recorded in Sloane's "The On-Line Encyclopedia of Integer Sequences (OEIS)" ("http://oeis.org/" accessed on March 8th, 2013) in Table (\ref{table}).

 Below we conclude the paper by listing some final properties that lead to new sequences whose proofs can be obtained by carefully implementing the definitions.

\begin{lemma}

\begin{enumerate}

\item    $N_{2\times 8}(j,k;j,1,1,1)=2^{(k-3)(j-1)}\cdot N_{2\times 8}(1,k;1,1,1,1) \text{ for } k\geq  3  \text{ and }  j\geq 1.$

\item  $  N_{2\times 8}(r,s;r,0,1,s-1)=N_{2\times 8}(r,s;r,0,s-1,1)=N_{2\times 8}(r,s;r,s-1,10)=N_{2\times 8}(r,s;r,1,s-1,0)=2^s-1 \text{ where } s\geq 2$.
For instance, $ N_{2\times 8}(r,3;r,0,2,1)=N_{2\times 8}(r,3;r,0,1,2)=N_{2\times 8}(r,3;r,1,2,0)=N_{2\times 8}(r,3;r,2,1,0)=7. $

\item $  N_{2\times 8}(r,s;r,0,k,s-k)= N_{2\times 8}(r,s;r,s-k,k,0).$  $  N_{2\times 8}(r,s;r,k,0,s-k)= N_{2\times 8}(r,s;r,s-k,0,k).$

\end{enumerate}

\end{lemma}

\begin{table}[h]
\begin{tabular}{|l|c|}\hline
The Sequences & Status \\\hline \hline
$N_{2\times 8}(\alpha,r;1,1,1,0);$ $ \qquad 1\leq \alpha \leq 8, 2\leq r \leq 4  $ &New \\ \hline
$N_{2\times 8}(r,2k;r,k,0,k)=\{6, 560, 714240, 13158776832, \ldots\}$  for $k\geq 1$  and $ r \geq 1  $ &New \\ \hline
$N_{2\times 8}(r+1,2;r,1,1,0)=\{36,84,180,372,756 ,...\};$ $ r \geq 1  $ &New \\ \hline
$N_{2\times 8}(r+1,3;r,1,1,1)=\{504,1176,2520,5208,10584,...\};$ $ r \geq 1  $ &New \\ \hline
$N_{2\times 8}(r+1,2r+1;r,0,r,r)=\{504,486080,1360627200,...\};$ $ r \geq 1  $ &New \\ \hline
$N_{2\times 8}(r+2,2r+1;r,0,1,r)=\{2352,9721600,449914060800,...\};$ $ r \geq 1  $ &New \\ \hline
$N_{2\times 8}(r,r+2;2,0,1,r)=\{840,52080,2187360,...\};$ $ r \geq 2  $ &New \\ \hline
$N_{2\times 8}(r,2k;r,k,k,0)=N_{2\times 8}(r,2k; r,0,k,k)=\{3,35,1395, 200787, \ldots\},  r,k\geq 1  $ & (*) \\ \hline
\end{tabular} \label{table} \caption{This table is partial list of the results. (*) exists in(\cite{numbers}) coded by A006098. }
\end{table}

\section{Conclusion}

In this work we established a formula that gives the number of distinct additive $\mathbb{Z}_{2} \mathbb{Z}_{8}$-codes. By specializing the parameters in this formula, we easily obtain the number of distinct codes over the ring $\mathbb{Z}_{8}$ and $\mathbb{Z}_{2} \mathbb{Z}_{4}$-codes. Further, some properties of this formula that is defined by the authors as Mixed Gaussian numbers are studied and some new number sequences are presented. Since Mixed Gaussian numbers are generalizations of Gaussian numbers we we believe that there so many properties some that are currently being investigated by the authors  that waits to be explored.


\begin{thebibliography}{9}

\bibitem{amis} I. Aydogdu and I. Siap, The Structure of $\mathbb{Z}_{2}\mathbb{Z}_{2^{s}}$ Additive Codes: Bounds on the minimum distance, \textsl{Applied Mathematics \& Information Sciences (AMIS)}, accepted (2013).

\bibitem{tilings}Jonathan Azose, A Tiling Interpretation of $q$-Binomial Coefficients,  \textsl{Hurvey Mudd College}, Department of Mathematics, Senior Thesis, (2007).

\bibitem{alger} Hacene Belbachir, Sadek Bouroubi, Abdelkader Khelladi, Connection between ordinary multinomials, Fibonacci numbers, Bell polynomials and discrete uniform distribution, \textsl{ Annales Mathematicae et Informaticae},   \textbf{35},   21–30, (2008).


\bibitem{deskiZ2Z4} Bilal, M., Borges, J., Dougherty, S., Fernandez, C., Optimal Codes over $\mathbb{Z}_{2}\mathbb{Z}_{4}$ \textsl{ In libro de acts  VII Jornadas de Matematica Discreta i Algoritmica}, Castro Urdiales (Spain),  131-139, (2010).



\bibitem{Optimal cod} Bilal, M., Borges, J., Dougherty, S., Fernandez, C., Extensions of $\mathbb{Z}_{2}\mathbb{Z}_{4}$ -additive self-dual codes preserving their properties, \textsl{IEEE International Symposium on Information Theory} , 3101	- 3105  (2012).


\bibitem{bona} Miklos Bona, Combinatorics of permutations,Discrete Mathematics and Its Applications, \textsl{Chapman and Hall/CRC}, (2004).

\bibitem{eskiZ2Z4} J. Borges, C. Fernandez, J. Pujol, J. Rifa and M. Villanueva, On $\mathbb{Z}_{2}\mathbb{Z}_{4}$ -linear codes
and duality, \textsl{ V Jornadas de Matematica Discreta i Algoritmica}, Soria (Spain), Jul. 11-14, 171-177, (2006).

\bibitem{Borges} J. Borges, C. Fern\'{a}ndez-C\'{o}rdoba, J. Pujol, J. Ri\'{f}a and M. Villanueva. $\mathbb{Z}_{2}\mathbb{Z}_{4}$-linear codes: generator matrices and duality, \textsl{Designs, Codes and
Cryptography}, \textbf{54} (2),  167-179, (2010).



\bibitem{mixed}Andries E. Brouwer, Heikki O. Hämäläinen, Patric R. J. Östergård, Neil J. A. Sloane,
Bounds on Mixed Binary/Ternary Codes, \textsl{ IEEE Transactions on Information Theory} \textbf{44} (1): 140-161 (1998)


\bibitem{Delsarte 1} Delsarte P., An algebraic approach to the association
schemes of coding theory, \textsl{Philips Research Rep.Supp.}, \textbf{10,}   vi+97   (1973).

\bibitem{Delsarte 2} Delsarte P,Levenshtein V.:Association schemes and
coding theory, \textsl{IEEE Trans.Inform.Theory},  \textbf{44}   (6)   2477-2504    (1998).

\bibitem{Hammons} Hammons AR, Kumar V, Calderbank AR, Sloane NJA, Sol\'{e}
P, The $%
\mathbb{Z}
_{4}$-linearity of Kerdock, Preparata, Goethals, and related codes, \textsl{IEEE
Trans. Inform. Theory},  \textbf{40}   301-319  (1994).


\bibitem{sloane} F.J. MacWilliams and N.J.A. Sloane, The Theory of
Error-Correcting Codes, North-Holland: New York, NY, (1977).



\bibitem{Vmulti} Amritanshu Prasad,  Counting Subspaces of a Finite Vector Space– 2,  \textsl{Resonance},  \textbf{15}, 12, (2010).


\bibitem{Pujol} Pujol J., Rif\'{a} J., Trranslation invariant propelinear
codes, \textsl{IEEE Trans. Inform. Theory}  \textbf{43 }590-598 (1997).

\bibitem{flag} C. Ryan Vinroot, An enumeration of ﬂags in ﬁnite vector spaces, \textsl{The Electronic Journal of Combinatorics},  \textbf{19}, 3, 1-9,  (2012).

\bibitem{S} Salturk E., Siap I., 	On Generalized Gaussian Numbers, \textsl{ Albanian Journal of Mathematics},   \textbf{6},  2  87-102 (2012).


\bibitem{Sa} Salturk E., Siap I., Generalized Gaussian Numbers Related to
Linear Codes over Galois Rings,\textsl{European Journal of Pure and Applied
Mathematics} \textbf{5}  250-259  (2012).

\bibitem{mk} Salturk E., Siap I., Generalized Gaussian Numbers and Some New Sequences, Physica Macedonica, accepted  \textbf{ 6}, (2013).

\bibitem{numbers} N. Sloane,   The On-Line Encyclopedia of Integer Sequences (OEIS), (\textsl{"http://oeis.org/"} accessed on March 8th, 2013).

\bibitem{ph}Hiroki Suyari, Mathematical structures derived from the q-multinomial coefficient in Tsallis statistics,  \textsl{Physica A: Statistical Mechanics and its Applications}, \textbf{  368}, 1, 63-82, (2006).






\end{thebibliography}
\end{document}